\documentclass[]{amsart}
 \usepackage[]{amsmath, amsthm, amsfonts,amssymb}
 \usepackage[all]{xy}
 \usepackage{hyperref}
 \xyoption{curve}

 \newcommand {\C} {{\mathbb C}}

 \newcommand {\Q} {{\mathbb Q}}
 \newcommand {\PP} {{\mathbb P}}

 \newcommand {\dt} {{\bullet}}

 \newcommand {\codim}{\text{codim}}
 \newcommand {\im}{\text{im}}

 \newtheorem {thm}[subsection]{Theorem}
 \newtheorem {cor}[subsection]{Corollary}
 \newtheorem {lemma}[subsection]{Lemma}

 \newtheorem {rmk}[subsection]{Remark}

\begin{document}
 \title{Hodge structure of a complete intersection of quadrics in a projective space}
 \author{Su-Jeong Kang}
 \address{Department of Mathematics\\
   Purdue University\\
   West Lafayette, IN 47907\\
   U.S.A.} 
\email{sjkang@math.purdue.edu}

 \maketitle

Let $V$ be a smooth projective variety of dimension $n$. One may say that $V$ has motivic dimension less than $d+1$ if the cohomology of $V$ comes from varieties of dimensions less than $d+1$ in some geometric way. More precisely, we say that $V$ has {\em motivic dimension} less than $d+1$ if there exists a (nonconnected) smooth projective variety $W$ of dimension less than $d+1$ and an algebraic correspondence $\Gamma$ on $W \times X$ such that $\Gamma$ induces a surjection $H^*(W) \to H^*(V)$.

Suppose the Generalized Hodge conjecture(GHC) (\cite{groth},\cite{lewis}) holds for $V$ and suppose the level of $V$ is less than $l$ (i.e. $\text{level}(H^i(V)) < l$ for all $i$). This implies that $V$ has a motivic dimension less than $l$. Conversely, if $V$ has motivic dimension less than $d+1$, i.e. there is a smooth projective variety $W$ of dimension less than $d+1$ and a surjection $\phi: H^*(W) \to H^n(V)$ induced by a correspondence $\Gamma$ on $W \times X$, then the level of the Hodge structure $H^n(V)$ is less than $d+1$ because a morphism of Hodge structures preserves the level. The existence of $W$ does not imply that the GHC holds for $V$ of course, because the dimension of the variety $W$ is not small enough to conclude the GHC for $V$. However, this gives a way to reduce the GHC for $V$ to the GHC for $W$, the variety with the smaller dimension. And we think that this can be a practical intermediate step for checking GHC for a smooth projective variety. 

For a smooth complete intersection of $k$ quadrics $V=Q_1 \cap \cdots \cap Q_k$ in $\PP^{n+k}$, it can be checked that the level of the Hodge structure $H^n(V)$ is less than $k$ (for example by checking Hodge number of it). In this paper, we show that a smooth complete intersection of $k$ quadrics in $\PP^{n+k}$ has motivic dimension less than $k $ (Theorem \ref{thm:main}). As a corollary of this, we get the Hodge conjecture for $V$ holds if $k<4$.

A brief outline of the proof of main theorem is as follows. By using the constuction in \cite{ogrady}, we form a family of quadrics parametrized by $\PP^{k-1}$ with the base locus $V$. In case when $n+k$ is odd, this is a family of even dimensional quadrics. Then we use the fact that an even dimensional quadric contains two irreducible families of linear spaces of the expected dimension, and we choose $W$ to be a double covering of $\PP^{n+k}$. In the case when $n+k$ is even, then quadrics in the original family are odd dimensional. So in order to find $W$, we pass to the family of singular fibers over the discriminant variety which can be understood as a family of quadrics of even dimension. Then we form a double covering of the discriminant variety to find $W$.

\medskip

All varieties in this paper will be defined over $\C$ and the cohomology without coefficient would be the singular cohomology with rational coefficient. 

\medskip

 I would like to express my thanks to Professor Donu Arapura for suggesting me to look at this problem and for his encouragement. The definition of motivic dimension is due to him. 


\section{Main Theorem}

Let $V =Q_1 \cap \cdots \cap Q_k$ be a smooth complete intersection of $k$ quadrics in a projective space $\PP^{n+k}$. Then $V$ is a smooth subvariety of $\PP^{n+k}$ of $\dim V =n$. Consider the family of quadrics $Q_t$ ($t \in \PP^{k-1}$) with the base locus $V$. We can give more precise description of this family. Let $Q_l = \{F_l=0 \}$ for $l=1,\cdots, k$, where $F_l(x_0, x_1,\cdots , x_{n+k}) = \sum^{n+k}_{i,j=0} c^l_{ij}x_i x_j$, $[c^l_{ij}]_{0 \leq i,j \leq n+k}$ is a symmetric $(n+k+1) \times (n+k+1)-$matrix. Then for a general $t=(t_0,..., t_{k-1}) \in  \PP^{k-1}$, the fiber $Q_t$ is given by the equation $\sum^{k-1}_{l=0}t_lF_{l+1}(x_0, ..., x_{n+k})=0$. Let $\Delta \subset \PP^{k-1}$ be the subvariety of $\PP^{k-1}$ parametrizing all singular fibers in the family. 
Then $\Delta$ is a hypersurface in $\PP^{k-1}$ of degree $n+k+1$. 
We assume that $\Delta$ is smooth. 
Set
$$
X = \{(t,x) \in \PP^{k-1} \times \PP^{n+k}~|~ x \in Q_t \subset \PP^{n+k} \} \subset \PP^{k-1} \times \PP^{n+k} 
$$
and let $p_1 : X \to \PP^{k-1}$ and $p_2 : X \to \PP^{n+k}$ be projections. Note that $X$ is a smooth projective variety of dimension $n+2k-2$.
\begin{align}\label{diag:proj}
\xymatrix@R=11pt{
X \ar[rr]^{p_2 \quad}  \ar[d]^{p_1}&& \PP^{n+k} \supset V \\ 
\PP^{k-1} && 
}
\end{align}

Let $i_X : X \hookrightarrow \PP^{k-1} \times \PP^{n+k}$ and $i_V: V \hookrightarrow \PP^{n+k}$ be inclusions. Then by Lefschetz theorem, the restriction maps $i^*_X : H^{n+2k-2}(\PP^{k-1} \times \PP^{n+k}) \to H^{n+2k-2}(X)$ and  $i^*_V:H^{n}(\PP^{n+k}) \to H^{n}(V)$ are injections. Set

$$
H^{n+2k-2}_0 (X) = H^{n+2k-2}(X) / \im ~ i^*_X 
$$
and 
$$
H^n_0 (V) = H^n (V) / \im ~i^*_V
$$

\begin{rmk}
In fact, 
$$
H^{n+2k-2}_0 (X) \cong \begin{cases}H^{n+2k-1}_c((\PP^{k-1}\times \PP^{n+k})-X) & \text{if $n$ is even}\\H^{n+2k-2}(X) & \text{if $n$ is odd}\end{cases}
$$
and 
$$
H^n_0(V) \cong \begin{cases}H^{n+1}_c(\PP^{n+k} -V) & \text{if $n$ is even}\\
H^n(V) & \text{if $n$ is odd}\end{cases}
$$
as Hodge structures.
\end{rmk}

Now we can state our main theorem precisely:

\begin{thm}\label{thm:main}
Let $V$ be a smooth complete intersection of $k$ quadrics in $\PP^{n+k}$. Then $V$ has motivic dimension less than $k$. More precisely, there is a smooth projective variety $W$ of dimension $k-1$ (resp. $k-2$) and surjection of rational Hodge structures
$$
\Theta : H^{k-1}(\tilde{W})(-q) \to H^n(V) \qquad \text{if $n$ is odd and $k$ is even}
$$
$$
(\text{resp. } \Theta : H^{k-2}(\tilde{W})(-q) \to H^n(V) \qquad \text{if $n$ is odd and $k$ is odd} )
$$
where $q=\frac{n-k+1}{2}$(resp. $\frac{n-k+2}{2}$) and $\tilde{W}$ is a disjoint union of finitely many copies of $W$, 
$$
\Theta : \bigoplus_r H^{2r}(W)(-q_r)^{\oplus l_r} \to  H^n_0(V) \qquad \text{if $n$ is even}
$$
where $q_r=\frac{n-2r}{2}$ and $l_r$ are positive integers given by explicit formula.
\end{thm}

We will prove this theorem by considering two cases depending on the parity of $n+k$ in the last two sections. One immediate corollary of this theorem is

\begin{cor}
If $k \leq 4$, then Hodge conjecture holds for $V$.
\end{cor}


\section{Intermediate step}
\label{sec:intermediate-step}

Consider the projection $p_2: X \to \PP^{n+k}$. Note that for any $q=(q_0,..., q_{n+k}) \in \PP^{n+k}- V$, 
\begin{eqnarray*}
  p^{-1}_2 (q) &=& \{ (t,q) \in \PP^{k-1} \times \PP^{n+k}~|~ q \in Q_t \}\\
&=& \left\{ t=(t_0,\cdots, t_{k-1})\in \PP^{k-1}~|~\sum^{k-1}_{l=0}\left( \sum^{n+k}_{i,j=0}  c^{l+1}_{ij} q_i q_j \right)t_l =0 \right\} \cong \PP^{k-2} 
\end{eqnarray*}
and  for any $q \in V$, 
$$
p^{-1}_2(q) = \{(t,q) \in \PP^{k-1} \times \PP^{n+k}~|~ q \in Q_t \} = \PP^{k-1} \times \{q\}
$$
since $V$ is the base locus of the family. Hence we have the following diagram, which we will use throughout this section:
\begin{align}\label{diag:xv}
\xymatrix{
(\PP^{n+k} - V) \times \PP^{k-2} \ar@{^(->}[r] \ar[d]^{p_2} & X \ar[d]^{p_2} &  \PP^{k-1} \times V \ar@{_(->}[l]_{i_E \quad} \ar[d]^{p_2} \\
\PP^{n+k} -V \ar@{^(->}[r] & \PP^{n+k} & \ar@{_(->}[l] V
}
\end{align}

 \begin{lemma}\label{lemma:rank1}
 $$
 \dim H^{n}_0(V) = \dim H^{n+2k-2}_0(X)
 $$
 \end{lemma}

\begin{proof}
From the diagram (\ref{diag:xv}), we get 
\begin{eqnarray}\label{eq:chiX}
  \chi(X) 
&=& (k-1)(n+k+1) + \chi(V)
\end{eqnarray}
Since $X$ is a very ample divisor in $\PP^{k-1} \times \PP^{n+k}$ and $V$ is a smooth complete intersection in $\PP^{n+k}$, by Lefschetz Theorem we have
$$
H^i(\PP^{k-1} \times \PP^{n+k}) \cong H^i(X), \qquad H^j(\PP^{n+k}) \cong H^j(V)
$$
for any $i< \dim X=n+2k-2$ and $j < \dim V=n$. 

\medskip

(Case 1) If $n$ is even :  then we have
\begin{eqnarray*}
  \chi(V)&=&2\sum^{n-1}_{i=0}(-1)^i b_i(V) + b_{n}(V)
=n + b_{n}(V)
\end{eqnarray*}
\begin{eqnarray*}  
\chi(X)&=&2\sum^{n+2k-3}_{i=0} (-1)^i b_i(X) + b_{n+2k-2}(X)\\
&=& 2 \left(  \frac{k(k+1)}{2} +k \cdot \frac{n-2}{2} \right) +b_{n+2k-2}(X)\\
&=& 
nk+k^2-k + b_{n+2k-2}(X)
\end{eqnarray*}
Therefore by (\ref{eq:chiX}) we get
$$
nk+k^2-k +b_{n+2k-2}(X)
=(k-1)(n+k+1)+n + b_{n}(V)\\
$$
i.e.
$$
b_{n+2k-2}(X) = b_{n}(V) + k-1
$$
 Note that 
$$\dim H^{n+2k-2}_0(X)=b_{n+2k-2}(X)-\dim H^{n+2k-2}(\PP^{k-1}\times \PP^{n+k})=b_{n+2k-2}(X)-k$$ 
and 
$$\dim H^{n}_0(V)=b_{n}(V)-\dim H^{n}(\PP^{n+k})=b_n(V) -1$$ 
therefore we get
 $$
 \dim H^{n+2k-2}_0(X)=\dim H^{n}_0(V)
 $$
in this case.
 
\medskip

(Case 2) If $n$ is odd : then by the similar calculation, we get 
$$
\chi(X) = nk+k^2 -b_{n+2k-2}(X), \qquad \chi(V) = (n+1)-b_{n}(V)
$$ 
Then, again by (\ref{eq:chiX})
\begin{eqnarray*}
  nk+k^2 - b_{n+2k-2}(X) 
= (k-1)(n+k+1)+(n+1)- b_{n}(V) =nk+k^2-b_{n}(V)
\end{eqnarray*}
i.e.
$$
b_{n+2k-2}(X) = b_n(V)
$$
Hence, we get
$$
\dim H^{n+2k-2}_0(X) = \dim H^{n+2k-2}(X) = \dim H^{n}(V)=\dim H^n_0(V)
$$
in this case also, which finishes the proof of the Lemma.
\end{proof}

\medskip

Set $E = \PP^{k-1} \times V \subset X$. Then we have a $\PP^{k-1}-$bundle $p_2: E \to V$. Let $i_E:E \hookrightarrow X$ be an inclusion. Then we have the following of morphism of Hodge structures 
$$
\phi :  H^{n}(V) \stackrel{p^*_2}{\longrightarrow} H^{n}(E) \stackrel{{i_E}_*}{\longrightarrow} H^{n+2k-2}(X)
$$
where ${i_E}_*: H^{n}(E) \to H^{n+2k-2}(X)$ is the Gysin map. (Note that $\codim(E,X)=k-1$)

\begin{thm}\label{thm:mainpartI}
$\phi: H^n(V) \to H^{n+2k-2}(X)$ induces a morphism
$$
\bar{\phi} : H^{n}_0 (V)(-k+1) \to H^{n+2k-2}_0(X)
$$
which is an isomorphism of rational Hodge structures.
\end{thm}

\begin{proof}
First we prove that the induced morphism 
$$
\bar{\phi}: H^{n}_0(V) \to H^{n+2k-2}_0(X)
$$
is well-defined. In the case when $n$ is odd, then $\bar{\phi}=\phi$, hence the morphism is well-defined. In case when $n$ is even, it is enough to show that $\phi(\im ~i^*_V) \subseteq \im~i^*_X$. Let $[E] \in H^{2k-2}(X)$ be the fundamental class of $E$. Note that $i^*_X: H^{2k-2}(\PP^{k-1}\times \PP^{n+k}) \to H^{2k-2}(X)$ is an isomorphism by Lefschetz theorem since $2k-2 < \dim X = n+2k-2$. Hence there is $\gamma \in H^{2k-2}(\PP^{k-1} \times \PP^{n+k})$ such that $i^*_X(\gamma) =[E]$. In fact, we can write $\gamma = [\PP^{k-1} \times S]$ where $[S] \in H^{2k-2}(\PP^{n+k})$ such that 
\begin{equation} \label{eq:gamma}
  i^*_X(\gamma) = [X \cap (\PP^{k-1} \times S)] = [\PP^{k-1} \times V]
\end{equation}
To show this, we consider the following commutative diagram:
$$
\xymatrix{
H^0(E) \ar[r]_{{i_*}_E} \ar@/^1pc/[rr]^{(i_X \circ i_E)_*} & H^{2k-2}(X) \ar[r]_{{i_X}_*} & H^{2k}(\PP^{k-1} \times \PP^{n+k}) \\
& H^{2k-2}(\PP^{k-1} \times \PP^{n+k}) \ar[u]^{i^*_X}_{\cong} \ar[ur]_{\cup [X]}^{\cong}
}
$$
Since $\gamma \in H^{2k-2}(\PP^{k-1} \times \PP^{n+k})$, we may write $\gamma=\sum^{k-1}_{i=0}[a_{i}H^i \times b_{i}L^{k-1-i}]$, where $H^i$(resp. $L^j$) is a linear space in $\PP^{k-1}$(resp. $\PP^{n+k}$) of codimension $i$(resp. $j$) By commutativity of the diagram, we have 
\begin{eqnarray*}
 [\PP^{k-1}\times V] ={i_X}_*[E]={i_X}_*{i^*}_X(\gamma) =\gamma \cup [X]=\sum^{k-1}_{i=0}[(a_{i} H^i \times b_{i}L^{k-1-i}) \cap X ]
\end{eqnarray*}
Hence $a_{i}=0$ for $i \neq 0$, since $H^i \subsetneqq \PP^{k-1}$ for $i \neq 0$ and then this forces $b_i =0$ for $i\neq 0$ since $[b_iL^{k-1-i}] \notin H^{2k-2}(\PP^{k-1} \times \PP^{n+k})$ for $i>0$. Hence   
$$
\gamma = [\PP^{k-1} \times b_0 L^{k-1}] 
$$
We take $S$ to be $b_0 L^{k-1}$ and we will use this $S$ later. 

Now consider the following diagram:
\begin{align}\label{diag:int}
\xymatrix@=10pt{
\ar@/^1pc/[rrrr]^{\phi} H^{n}(V)\ar@{}[ddrr]|{(I)} \ar[rr]^{p^*_2} && H^n(E) \ar@{}[drr]|{(III)} \ar[rr]^{{i_E}_*} && H^{n+2k-2}(X) \\
&&&  H^n(X) \ar[ul]_{i^*_E} \ar[ur]^{\cup [E]} & \\
H^{n}(\PP^{n+k}) \ar[uu]^{i^*_V} \ar[rr]_{{pr}^*_2} && H^{n}(\PP^{k-1} \times \PP^{n+k}) \ar@{}[uur]|{(II)\qquad }   \ar[ur]_{i^*_X} \ar[uu]^{j^*_E} \ar[rr]_{\Phi} && H^{n+2k-2}(\PP^{k-1} \times \PP^{n+k}) \ar[uu]_{i^*_X} \ar@{}[uul]|{\qquad (IV)}
}
\end{align}
where
\begin{enumerate}
\item  $j_E: E \to \PP^{k-1} \times \PP^{n+k}$,~ $i_E: E \to X$ and $i_X:X \to \PP^{k-1} \times \PP^{n+k}$ are inclusions,
\item  $pr_2: \PP^{k-1}\times \PP^{n+k} \to \PP^{n+k}$ is the projection to the second factor,
\item $\Phi:H^n(\PP^{k-1} \times \PP^{n+k}) \to H^{n+2k-2}(\PP^{k-1} \times \PP^{n+k})$ is defined by $\Phi(\alpha) = \alpha \cup \gamma$.
\end{enumerate}
Commutativity of $(II)$ and $(III)$ are clear and commutativity of $(I)$ follows from the following commutative diagram:
$$
\xymatrix{
E = V \times \PP^{n+k} \ar@{^(->}[rr]^{\quad j_E} \ar[d]^{p_2} && \PP^{k-1} \times \PP^{n+k} \ar[d]^{pr_2} \\
V \ar@{^(->}[rr]^{i_V} && \PP^{n+k}
}
$$
We show the commutativity of $(IV)$. Let $\alpha \in H^{n}(\PP^{k-1} \times \PP^{n+k})$. Then 
$$
i^*_X \circ \Phi (\alpha)= i^*_X(\alpha \cup \gamma)=i^*_X(\alpha) \cup i^*_X(\gamma) =i^*_X(\alpha) \cup [E]
$$
by the definition of $\gamma$. Hence $(IV)$ commutes. Therefore we have
$$
 \phi \circ i^*_V = i^*_X \circ \Phi \circ pr^*_2 
$$
and hence $\phi(\im ~ i^*_V) \subseteq  \im ~ i^*_X$ and $\phi$ induces a well-defined morphism $\bar{\phi}: H^n_0(V) \to H^{n+2k-2}_0(X)$ in this case also.

\medskip

To show that $\bar{\phi}$ is an isomorphism, note that we have
\begin{eqnarray}\label{eq:E}
H^{n}(E)&\cong& ( H^{n}(V) \otimes H^0(\PP^{k-1}) ) \oplus \left(\bigoplus^{k-1}_{i=1} H^{n-2i} (\PP^{n+k}) \otimes H^{2i}(\PP^{k-1}) \right) 
\end{eqnarray}
by the K\"{u}nneth formula and Lefschetz theorem.

\medskip
(Case 1) If $n$ is odd : in this case $\phi=\bar{\phi}$ and we show that $\phi={i_E}_* \circ p^*_2$ is surjective. Since $n$ is odd, $H^{n-2i}(\PP^{k-1}) =0$ for all $i $. Hence (\ref{eq:E}) gives an isomorphism
$$
p^*_2 : H^{n}(V)  \stackrel{\cong}{\longrightarrow} H^{n}(E)
$$
Now for the morphism ${i_E}_*: H^{n}(E) \to H^{n+2k-2}(X)$, consider the following Gysin exact sequence
\begin{equation}\label{eq:gysin}
\cdots \to H^{n+2k-3}(X-E) \to H^{n}(E) \stackrel{{i_E}_*}{\to} H^{n+2k-2}(X) \to H^{n+2k-2}(X-E) \to\cdots
\end{equation}
Since $X -E \cong \PP^{k-2} \times (\PP^{n+k} -V)$,
\begin{eqnarray*}
H^{n+2k-2}(X-E)&=&H^{n+2k-2} (\PP^{k-2} \times (\PP^{n+k} -V))\\ 
&=& \bigoplus^{k-2}_{i=0} H^{(n+2k-2)-2i}(\PP^{n+k}-V)(-i)=0
\end{eqnarray*}
Hence, ${i_E}_* : H^{n}(E) \to H^{n+2k-2}(X)$ is surjective and hence we get the surjection $\phi = {i_E}_* \circ p^*_2$. Now lemma \ref{lemma:rank1} implies that $\bar{\phi}=\phi$ is an isomorphism in this case. 

\medskip

(Case 2) If $n=2l$ is even : then (\ref{eq:E}) gives that $p^*_2: H^{n}(V) \to H^{n}(E)$ is an injection. 
Since $X -E \cong \PP^{k-2} \times (\PP^{n+k} -V)$, we get
$$
H^{n+2k-3}(E)=H^{n+2k-3} (\PP^{k-2} \times (\PP^{n+k} -V)) = \bigoplus^{k-2}_{i=0} H^{(n+2k-3)-2i}(\PP^{n+k}-V)(-i)=0
$$
Then Gysin exact sequence (\ref{eq:gysin}) implies that ${i_E}_* : H^{n}(E) \to H^{n+2k-2}(X)$ is injective. Since $\phi= {i_E}_* \circ p^*_2$, we get the injection $\phi = {i_E}_* \circ p^*_2 : H^{n}(V) \to H^{n+2k-2}(X)$. We show that the induced map $\bar{\phi}$ is also an injection in this case: Suppose $\bar{\phi}(\bar{\alpha})=0$. Then 
$$
\phi(\alpha)\in \im[i^*_X: H^{n+2k-2}(\PP^{k-1}\times \PP^{n+k})\to H^{n+2k-2}(X)]$$
where $\alpha \in H^{n}(V)$ which maps to $\bar{\alpha} \in H^{n}_0(V)$. Let $\beta \in H^{n+2k-2}(\PP^{k-1}\times \PP^{n+k})$ such that $i^*_X(\beta) = \phi(\alpha)$. We claim that $\beta$ is a cycle in $H^0(\PP^{k-1}) \otimes H^{n+2k-2}(\PP^{n+k})$.
$$
\xymatrix@=12pt{
 &&& H^{n+2k-2}(E)\\
H^n(V) \ar[r]^{p^*_2}&H^n(E)  \ar[r]^{{i_E}_*} \ar@/^1pc/[urr]^{\cup c_{k-1}(N_{E/X})} & H^{n+2k-2}(X) \ar[ur]^{i^*_E}& \\
H^n(\PP^{n+k}) \ar[u]^{i^*_V} \ar[r]_{pr^*_2}& \ar[u]H^{n}(\PP^{k-1}\times \PP^{n+k})\ar[r]_{\Phi}& H^{n+2k-2}(\PP^{k-1}\times \PP^{n+k}) \ar[u]^{i^*_X} \ar[uur]_{(i_X \circ i_E)^*}&
}
$$
By applying $i^*_E$ to $\phi(\alpha) =i^*_X(\beta) \in H^{n+2k-2}(X)$, we get
\begin{multline}\label{eq:phialpha}
i^*_E(\phi(\alpha))= i^*_E {i_E}_* p^*_2(\alpha) = p^*_2(\alpha) \cup c_{k-1}(N_{E/X}) 
= i^*_E(i^*_X(\beta) )) =(i_X \circ i_E)^*(\beta)
\end{multline}
where the second equality comes from the self-intersection formula (\cite[p103]{fulton}). From the inclusions $E \stackrel{i_E}{\hookrightarrow} X \stackrel{i_X}{\hookrightarrow}\PP^{k-1} \times \PP^{n+k}$, we have
$$
0 \to i^*_E (N_{X/\PP^{k-1}\times \PP^{n+k}}) \to N_{E/\PP^{k-1}\times \PP^{n+k}} \to N_{E/X} \to 0
$$
Hence 
$$
c ( N_{E/\PP^{k-1}\times \PP^{n+k}}) = c( i^*_E (N_{X/\PP^{k-1}\times \PP^{n+k}})) \cdot c(N_{E/X})
$$
In particular,
\begin{multline}\label{eq:ck0}
c_{k} ( N_{E/\PP^{k-1}\times \PP^{n+k}}) \\= c_1( i^*_E (N_{X/\PP^{k-1}\times \PP^{n+k}})) \cdot c_{k-1}(N_{E/X})= i^*_E([X]|_X) \cup_E [E]|_E
\end{multline}
since $X$ is a divisor in $\PP^{k-1} \times \PP^{n+k}$, where $\cup_E$ is the cup product on $E$. Since the inclusion $i_X \circ i_E: E \to \PP^{k-1} \times \PP^{n+k}$ is actually $(\text{id}_{\PP^{k-1}}, i_V) :\PP^{k-1} \times V \hookrightarrow \PP^{k-1}\times \PP^{n+k}$, we have $N_{E/\PP^{k-1}\times \PP^{n+k}}=p^*_2(N_{V/\PP^{n+k}})$ and hence
\begin{equation}\label{eq:ck}
c_{k} ( N_{E/\PP^{k-1}\times \PP^{n+k}}) =p^*_2(c_k(N_{V/\PP^{n+k}})) \in H^0(\PP^{k-1}) \otimes H^{2k}(V) \subset H^{2k}(E)
\end{equation}
Now since $X$ is an ample divisor in $\PP^{k-1} \times \PP^{n+k}$, $[X] = [a H + b L] \in H^2(\PP^{k-1}\times \PP^{n+k})$ where $H$(resp. $L$) is a hyperplane in $\PP^{k-1}$(resp. $\PP^{n+k}$) and $a,b \geq 0$ such that $ab\neq 0$. Now by using (\ref{eq:gamma}),
\begin{eqnarray*}
 i^*_E([X]|_X) \cup_E [E]|_E &=& i^*_E ([X]|_X \cup_X [E])=i^*_E ([X]|_X \cup_X i^*_X (\gamma))\\
&=&i^*_E i^*_X([X] \cup_{\PP^{k-1} \times \PP^{n+k}} \gamma)\\
&=&  i^*_E i^*_X ([a H +b L ]\cup_{\PP^{k-1} \times \PP^{n+k}} [\PP^{k-1} \times S])\\
&=& i^*_E i^*_X ([a (H \times S) + b (\PP^{k-1} \times (S \cap L))])\\
 &\in & (H^2(\PP^{k-1})\otimes H^{2k-2}(V)) \oplus (H^0(\PP^{k-1}) \otimes H^{2k}(V))
\end{eqnarray*} 
Then by (\ref{eq:ck0}), (\ref{eq:ck}), we have $a=0$ and hence $[X] = [bL]\in H^2(\PP^{k-1}\times \PP^{n+k})$. Set $P$ be a hypersurface of degree $b$ in $\PP^{n+k}$ such that $[X] =[\PP^{k-1}\times P]$. 
Since $X$ contains $E=\PP^{k-1} \times V$, we may assume that $V \subset P$. Then $N_{E/X} = p^*_2 (N_{V/P})$, and hence 
\begin{equation*}
c_{k-1}(N_{E/X}) = p^*_2(c_{k-1}(N_{V/P})) \in H^0(\PP^{k-1})\otimes H^{2k-2}(V)
\end{equation*}
Now, since $p^*_2(\alpha) \in H^0(\PP^{k-1}) \otimes H^{n}(V)$, we have
$$
i^*_E \phi (\alpha) =  p^*_2(\alpha)  \cup c_{k-1}(N_{E/X})  \in H^0(\PP^{k-1})\otimes H^{n+2k-2}(V)
$$
Then by (\ref{eq:phialpha}), we have
$$
i^*_E \phi (\alpha)= (i_E \circ i_X)^*(\beta) = (\text{id}_{\PP^{k-1}}, i^*_V)(\beta) \in H^{0} (\PP^{k-1})\otimes H^{n+2k-2}(V)
$$
Hence $\beta \in H^0(\PP^{k-1}) \otimes H^{n+2k-2}(\PP^{n+k})$.

Next note that by cupping with $\gamma = [\PP^{k-1} \times S]$, $\Phi: H^{n}(\PP^{k-1}\times \PP^{n+k}) \to H^{n+2k-2}(\PP^{k-1}\times \PP^{n+k})$ maps the K\"{u}nneth component $H^{i}(\PP^{k-1})\otimes H^{n-i}(\PP^{n+k})$ of $H^n(\PP^{k-1}\times \PP^{n+k})$ to the K\"{u}nneth component $H^{i}(\PP^{k-1}) \otimes H^{n+2k-2-i}(\PP^{n+k})$ of $H^{n+2k-2}(\PP^{k-1} \times \PP^{n+k})$ isomorphically. Hence we have $\eta \in H^0(\PP^{k-1}) \otimes H^{n}(\PP^{n+k})$ such that $\eta \cup \gamma = \Phi(\eta)= \beta$. Now $\eta \in \im [pr^*_2: H^n(\PP^{n+k}) \to H^{n}(\PP^{k-1}\times \PP^{n+k})]$ and we may consider $\eta \in H^{n}(\PP^{n+k})$ since $pr^*_2$ is injective. Then by commutativity of diagram (\ref{diag:int}), we have
$$
\phi \circ i^*_V(\eta)= i^*_X \circ \Phi \circ {pr^*_2}(\eta)= i^*_X(\eta \cup \gamma)=i^*_X(\beta) =\phi(\alpha)
$$
Since $\phi$ is injective, we have $\alpha = i^*_V(\eta)$. Hence $\bar{\phi}$ is also an injection. Now by lemma \ref{lemma:rank1}, we can conclude that $\bar{\phi}: H^{n}_0(V) \to H^{n+2k-2}_0(X)$ is an isomorphism.
\end{proof}

\section{Proof of Theorem \ref{thm:main} when $n+k$ is odd}
\label{sec:part-ii-n}

Throughout this section, we assume that $n+k=2m+1$.

\medskip
In order to prove theorem \ref{thm:main}, we use the construction of O'Grady\cite{ogrady}. We give a brief outline of his construction here. For detailed construction, see \cite{ogrady}.

\medskip

Recall the diagram (\ref{diag:proj}) and consider the projection $p_1 : X \to \PP^{k-1}$ and recall that $\Delta$ is the discriminant variety, which is a smooth hypersurface in $\PP^{k-1}$ by our assumption. In case when $n+k$ is odd, for a general $t \in \PP^{k-1}$, the fiber $p^{-1}_1(t) = Q_t$ is a smooth quadric of dimension $n+k-1=2m$ in $\PP^{n+k}$. Hence it contains two irreducible families of $m-$planes parametrized by $F^1_t$ and $F^2_t$. Note that $F^1_t \cong F^2_t$ and $\dim F^i_t = \frac{m(m+1)}{2}$ \cite{gh}. Let $F$ be the abstract variety to which $F^i_t$ is isomorphic for $i=1,2$ and for any $t$. Let $W$ be a double covering of $\PP^{k-1}$ branched over $\Delta$ and let $\sigma: W \to \PP^{k-1}$ be the covering map. Set
$$
P =\{ M \subset X ~|~ p_1(M) = t \in \PP^{k-1} \text{ a point, } p_2(M) \cong \PP^{m} \subset \PP^{n+k}\}
$$ 
Then there is a natural map $\psi : P \to \PP^{k-1}$ defined by $\psi(M) = p_1(M) \in \PP^{k-1}$. Then the Stein factorization of $\psi$ is factored through $W$ and get a composition
$$
\psi  : P \stackrel{f}{\longrightarrow} W \stackrel{\sigma}{\longrightarrow} \PP^{k-1}
$$
and set
\begin{eqnarray*}
\Gamma &=& \{ (M,x) \in P \times X ~|~ x \in M \subset X\} \\ 
&=& \{(M,x) \in P \times X ~|~ p_1(x) = p_1(M) \in \PP^{k-1},~p_2(x) \in p_2(M)\cong \PP^{m} \} 
\end{eqnarray*}
Let $pr_1 : \Gamma \to P$ and $pr_2 : \Gamma \to X$ be the projections. We summarize the construction in the following diagram:
$$
\xymatrix@R=10pt{
& \Gamma \ar[dl]_{pr_1} \ar[dr]^{pr_2}& \\
P \ar[d]^f \ar[drr]^{\psi} && X \ar[d]^{p_1} \\
W \ar[rr]_{\sigma} && \PP^{k-1}
}
$$
Note that for any $w \in W$,
$$
f^{-1}(w) = \{ M \in X~|~ p_1(M) = \sigma(w), p_2(M) \cong \PP^{m} \in \PP^{n+k}\}
\cong\{ \PP^m \subset Q_{\sigma(w)} \} \cong F
$$
and for any $M \in P$, 
\begin{eqnarray*}
pr^{-1}_1(M) &=& \{(M,x) \in P \times X~|~ x \in M \subset X \}\\
&\cong& \{x \in X ~|~ p_1(x) =p_1(M),~ p_2(x) \in p_2(M) \cong \PP^{m} \} \cong \PP^{m}
\end{eqnarray*}
Hence $f:P \to W$ is $F-$bundle and $pr_1: \Gamma \to P$ is $\PP^m-$bundle. 

 \medskip

 Now we give a proof of theorem \ref{thm:main} in case of $n+k$ odd. 

\begin{proof}[Proof of theorem \ref{thm:main} when $n+k$ is odd]

Consider the morphism of Hodge structures
$$
{pr_2}_* : H_{n+2k-2}(\Gamma) \to H_{n+2k-2}(X)
$$
We show that ${pr_2}_*$ is a surjection. First note that $pr_2:\Gamma \to X$ is a surjection. Hence we can take an iterated hyperplane section $\Gamma_1$ of $\Gamma$ such that $\dim \Gamma_1 = \dim X = n+2k-2$ and $\Gamma_1$ surjects onto $ X$. 
Let $g = pr_2|_{\Gamma_1}:\Gamma_1 \to X$. Then $g$ is a generically finite map. Let $j : \Gamma_1 \to \Gamma$ be an inclusion. Then we have the following commutative diagram:

$$
\xymatrix{
\Gamma \ar[dr]^{pr_2}& &&H_{n+2k-2}(\Gamma) \ar[drr]^{{pr_2}_*} && \\
\Gamma_1 \ar@{^(->}[u]^j \ar[r]^{g} & X && H_{n+2k-2}(\Gamma_1) \ar[u]^{j_*} \ar[rr]^{g_*} \ar[d]^{P_{\Gamma_1}}_{\cong} && H_{n+2k-2}(X)  \ar[d]^{P_{X}}_{\cong}\\
 &&&H^{n+2k-2}(\Gamma_1) \ar[rr]^{g_*} && H^{n+2k-2}(X) \ar@/^1pc/[ll]^{g^*} 
}
$$ 
where $P_{\Gamma_1}$ and $P_X$ are isomorphisms from Poincar\'{e} duality. Since $g_* ={pr_2}_* \circ j_*$, it is enough to show that $g_*$ is surjective, but it is clear since $g_* \circ g^* = \deg g \cdot \text{id}$. 

Now since $pr_1: \Gamma \to P$ is $\PP^m-$bundle, by K\"{u}nneth formula we have
\begin{eqnarray*}
H_{n+2k-2}(\Gamma) &=& \bigoplus^{m}_{j=0} ( H_{n+2k-2-2j}(P) \otimes H_{2j}(\PP^m)) 
\\
&=& \bigoplus^{m}_{j=0} H_{n+2k-2-2j}(P)(j)
\end{eqnarray*}
We claim that for each $j$, there is an isomorphism 
$$
H_{n+2k-2-2j}(P) = \bigoplus_{r+2s=n+2k-2-2j}(H_{r}(W) \otimes H_{2s}(F)) 
$$ 
To show this, first note that we can use cohomology instead of homology since all our varieties considered are smooth. Let $U_W \subset W$ be a Zariski open set in $W$ such that $f^{-1}(U_W) \cong U_W \times F$. Since $F$ has a celluar decomposition \cite{ehr}, $H_*(F)$ are generated by algebraic cycles. Let $\alpha_1 , \cdots, \alpha_l$ be the algebraic cycles generating $H_*(F)$. i.e. there are algebraic subvarieties $Z_1, \cdots, Z_l$ of $F$ such that the fundamental classes of them are $\alpha_1,..., \alpha_l$. Let $ p: U_W \times F \to F$ be the projection to the second factor and consider the algebraic cycles $p^*(\alpha_i)$ in $H^*(U_W \times F)$ which are supported on $U_W \times Z_i$ for $i=1,2,...,l$. Let $\beta_i$ be the closures of $p^*(\alpha_i)$ in $P$ for each $i=1,..., l$. Then this gives the splitting of the restriction map $H^*(P) \to H^*(F)$. So we may apply the Leray-Hirsch theorem \cite{sp}. Therefore we have
$$
H^q(P) = \bigoplus_{r+2s=q}(H^r(W) \otimes H^{2s}(F))
$$
for any $q$. In particular,
\begin{equation}\label{eq:p}
H^{2d-(n+2k-2-2j)}(P) = \bigoplus_{s}(H^{2d-(n+2k-2-2j)-2s}(W) \otimes H^{2s}(F))
\end{equation}
where $d =\dim P$.

\medskip

(Case 1) If $n$ is odd (and hence $k$ is even) : in this case, from (\ref{eq:p}) we have
$$
H^{2d -(n+2k-2-2j)}(P ) \cong H^{k-1}(W)\otimes H^{2d-(n+3k-3-2j)}(F)
$$
for each $j$, since $W$ is simply connected, $H^r(W)=0$ for all odd $r$ such that $r \neq \dim W =k-1$, or equivalently
$$
H_{n+2k-2-2j}(P) \cong H_{k-1}(W) \otimes H_{n+k-1-2j}(F) = H_{k-1}(W)(m-j)^{\oplus l_j}
$$
where $l_j = \dim H_{n+k-1-2j}(F)$. Therefore we have
\begin{eqnarray*}
  H_{n+2k-2}(\Gamma) &=& \bigoplus^{m}_{j=0} H_{n+2k-2-2j}(P)(j)\\
&\cong& \bigoplus^{m}_{j=0} \left(H_{k-1}(W)(m-j )^{\oplus l_j} \right)(j) = H_{k-1}(W)(m)^{\oplus N}
\end{eqnarray*}
where $N = \sum^m_{j=0} l_j$. Thus, we have a surjection
$$
{pr_2}_* : H_{k-1}(W)(m)^{\oplus N} \to H_{n+2k-2}(X)
$$
and by Poincar\'{e} duality, we get a surjection
$$
\xymatrix{
H^{k-1}(W)(-m)^{\oplus N} \ar@{>>}[r] & H^{n+2k-2}(X) 
}
$$
By choosing $\tilde{W}$ to be the disjoint union of $N$ copies of $W$ and by composing with $\bar{\phi}^{-1}$ in theorem \ref{thm:mainpartI}, we get a surjection
$$
\xymatrix{
\Theta : H^{k-1}(\tilde{W})(-m+k-1) \ar@{>>}[r] & H^n (V)
}
$$
as we claimed.

\medskip

(Case 2) If $n$ is even (and hence $k$ is odd) : in this case, (\ref{eq:p}) gives 
$$
H^{2d -(n+2k-2-2j)}(P ) \cong \bigoplus_{r} (H^{2r}(W) \otimes H^{2d-(n+2k-2-2j)-2r}(F))
$$
for each $j$, or equivalently
$$
H_{n+2k-2-2j}(P) \cong \bigoplus_s \left(H_{2s}(W) \otimes H_{n+2k-2-2j-2s}(F) \right)
$$
 Therefore we have
\begin{eqnarray*}
  H_{n+2k-2}(\Gamma) &=& \bigoplus^{m}_{j=0} H_{n+2k-2-2j}(P)(j)\\
&\cong&  \bigoplus^{m}_{j=0}\bigoplus_s \left(H_{2s}(W) \otimes H_{n+2k-2-2j-2s}(F) \right)(j)\\
&=& \bigoplus_s H_{2s}(W) \otimes \left( \bigoplus^m_{j=0} H_{n+2k-2-2j-2s}(F) \right)(j)\\
&=& \bigoplus_s H_{2s}(W)(q_s)^{\oplus l_s}
\end{eqnarray*}
where $q_s = \frac{n+2(k-1)-2s}{2}$ and $l_s = \sum^m_{j=0} \dim H_{n+2k-2-2j-2s}(F)$. Thus, we have a surjection
$$
\xymatrix{
{pr_2}_* : \bigoplus_s H_{2s}(W)(q_s)^{\oplus l_s} \ar@{>>}[r] & H_{n+2k-2}(X)
}
$$
Recall that $H^{n+2k-2}_0(X) = H^{n+2k-2}(X)/ \im ~i^*_X$. Hence we get a composition of surjections
$$
\xymatrix{
\bigoplus_s H^{2s}(W)(-q_s)^{\oplus l_s} \ar@{>>}[r] & H^{n+2k-2}(X) \ar@{>>}[r] & H^{n+2k-2}_0(X)
}
$$
Therefore, by composing with $\bar{\phi}^{-1}$ from theorem \ref{thm:mainpartI} we get a surjection of Hodge structures
$$
\xymatrix{
\Theta : \bigoplus_s H^{2s}(W)(-q_s+k-1)^{\oplus l_s} \ar@{>>}[r] & H^{n}_0(V)
}
$$
(Note that $q_s -k+1 = \frac{n-2s}{2}$.) This completes the proof of theorem \ref{thm:main} in case when $n+k$ is odd.
\end{proof}

\section{Proof of Theorem \ref{thm:main} when  $n+k$ is even}
\label{sec:case-when-n}

Throughout this section, we assume that $n+k=2m$ is even.  

\medskip

Again we start by considering the projection $p_1: X \to  \PP^{k-1}$ in diagram (\ref{diag:proj}). Recall that the discriminant variety $\Delta$ is a smooth hypersurface of degree $n+k+1$ in $\PP^{k-1}$ by our assumption. Set $U =\PP^{k-1} -\Delta$ and let $X_{\Delta} = p^{-1}_1(\Delta)$ and $X_{U} =p^{-1}_1(U)$.
\begin{align}\label{diag}
\xymatrix{
X_{\Delta} \ar@{^(->}[r] \ar[d]^{p_1} & X \ar[d]^{p_1} & \ar@{_(->}[l] \ar[d]^{p_1} X_U \\
\Delta \ar@{^(->}[r] & \PP^{k-1} &\ar@{_(->}[l] U =\PP^{k-1} - \Delta
}
\end{align}
Note that for any $t \in \Delta$, $p^{-1}_1(t)$ is a singular quadric in the family. Since we have assumed that $\Delta$ is smooth, a singular fiber is a cone through a $0-$plane (i.e. a point) over a quadric of rank $n+k$ in $\PP^{n+k-1}$, i.e. all singular fibers are cones over a smooth quadric $\tilde{Q}_t$ of dimension $n+k-2$ in $\PP^{n+k-1}$. For any $t\in \Delta$, we denote $p^{-1}_1(t) = C_t$ a cone with a vertex $0_t$. We can form a family of quadrics of dimension $n+k-2$ over $\Delta$ as follows: Let $s: \Delta \to X_{\Delta}$ be a section of $p_1$ defined by $s(t) = 0_t$ for any $t\in \Delta$. Let $Y_0$ be a general hyperplane section of $X_{\Delta} - s(\Delta)$ and let $p_1|_{Y_0}: Y_0 \to \Delta$ be the obvious map. Let $Y$ be a smooth compactification of $Y_0$. Then by using theorem by Hironaka, we may assume that the rational map $Y \to \Delta$ is actually a morphism. We denote this morphism by $\pi:Y \to \Delta$. Then for a general $t \in \Delta$, $\pi^{-1}(t) = \tilde{Q}_t$ a smooth quadric of dimension $n+k-2$ in $\PP^{n+k-1}$.
$$
\xymatrix{
Y   \ar[dr]_{\pi} & \ar@{_(->}[l]Y_0 \ar@{^(->}[r] \ar[d]^{p_1|_{Y_0}}  & X_{\Delta} \ar[d]^{p_1}  \ar@{^(->}[r] & X \ar[d]_{p_1}\\
&\Delta \ar@{=}[r]&\Delta \ar@{^(->}[r] & \PP^{k-1}
}
$$

\begin{lemma}
For any $p$
  $$
H^{p}(X_{\Delta}) \cong H^{p-2}(Y)(-1) \quad \text{as Hodge structures}
$$
\end{lemma}

\begin{proof}
Consider the Leray spectral sequence associated to the map $p_1: X_{\Delta} \to \Delta$
$$
'E^{pq}_2 = H^p(\Delta, R^q {p_1}_* \Q ) \quad \Rightarrow \quad H^{p+q}(X_{\Delta},\Q)
$$
and one associated to the map $\pi : Y \to \Delta$
$$
''E^{pq}_2 = H^p (\Delta, R^q {\pi}_* \Q) \quad \Rightarrow \quad H^{p+q}(Y,\Q)
$$
Since for any $t \in \Delta$, $C_t$ is a cone through a point $0_t$ over a smooth quadric $\tilde{Q}_t$ in $\PP^{n+k-1}$, we have 
\begin{eqnarray*}
R^q {p_1}_* \Q &=& R^q {(p_1|_{Y_0})}_! \Q = (R^{2(n+k-1)-q}{({p_1|}_{Y_0})}_*\Q)^*(-n-k+1) \\
&=& (R^{2(n+k-1)-q}{\pi}_*\Q)^*(-n-k+1) = R^{q-2}{\pi}_*\Q(-1)
\end{eqnarray*}
Hence we have
$$
'E^{pq}_2 ~=~ ''E^{p,q-2}_2(-1)
$$
and so
$$
H^{p+q}(X_{\Delta}) \cong H^{p+q-2}(Y)(-1)
$$
as Hodge sturctures \cite{arapura}. 
\end{proof}

In particular,
\begin{equation}\label{eq:9}
H^{n+2k-2}(X_{\Delta}) \cong H^{n+2k-4}(Y)(-1)
\end{equation}

\begin{lemma}\label{lemma:inj}
  There is an injection of Hodge structures 
$$
0 \to H^{n+2k-2}_0 (X) \to H^{n+2k-4}(Y)(-1) 
$$
\end{lemma}

 \begin{proof}

From the top row of the diagram (\ref{diag}), we have an exact sequence of mixed Hodge structures

\begin{equation}\label{eq:mhs}
\cdots \to H^{n+2k-2}_c (X_U) \to H^{n+2k-2}(X) \to H^{n+2k-2}(X_{\Delta}) 
\to \cdots
\end{equation}

By using (\ref{eq:9}) and by taking the exact functor $Gr^W_{n+2k-2}$, we get

\begin{equation}
  \label{eq:7}
  0 \to Gr^W_{n+2k-2}  H^{n+2k-2}_c (X_U) \to H^{n+2k-2}(X) \to H^{n+2k-4}(Y)(-1) \to \cdots
\end{equation}

  Consider the morphism $p_1: X_U \to U$ and the Leray spectral sequence associated to it
$$
E^{pq}_2 = H^p_c (U , R^q {p_1}_* \Q) \quad \Rightarrow \quad H^{p+q}_c(X_U,\Q)
$$
Note that
$$
E^{pq}_2 = \begin{cases} 0 & \text{if $q$ is odd} \\ H^p_c(U,\Q) & \text{if $q$ is even} \end{cases}
$$
since $(R^q {p_1}_* \Q)_t = H^q (Q_t,\Q)$ and $Q_t$ is a smooth quadric of dimension $n+k-1$(odd). Also, we have 
$$
0 \to H^0_c(U) \to H^0(\PP^{k-1}) \to H^0(\Delta) \to H^1_c(U) \to 0
$$
$$
0 \to H^{2i-1}(\Delta) \to H^{2i}_c(U) \to H^{2i}(\PP^{k-1}) \to H^{2i}(\Delta) \to H^{2i+1}_c(U) \to 0  
$$
for $1 \leq i \leq k-2$ and 
$$
H^{2k-2}_c(U) \cong H^{2k-2}(\PP^{k-1}) 
$$
Now, since $\Delta$ is a smooth hypersurface in $\PP^{k-1}$, by Lefschetz theorem we get
$$
H^j(\PP^{k-1}) \stackrel{\cong}{\longrightarrow} H^j(\Delta), \qquad \text{for } j < \dim \Delta= k-2
$$
Hence $H^{2i-1}(\Delta)  = 0$ and hence $H^{2i}_c(U)=H^{2i+1}_c(U)=0$ for $i$ such that $2i < \dim \Delta = k-2$. By applying duality on $H^{2i-1}(\Delta)$, we get $H^j_c(U)=0$ unless $j = k-1$ or $j=2k-2$. Hence the Leray spectral sequence degenerates at $E_2$ and we have a short exact sequence
$$
0 \to E^{2k-2,n}_{\infty} \to H^{n+2k-2}_c(X_U) \to E^{k-1,n+k-1}_{\infty} \to 0
$$
Note that $E^{k-1,n+k-1}_{\infty} = E^{k-1,n+k-1}_2 =0$ since $n+k$ is even. Hence we have
\begin{equation}\label{eq:short}
E^{2k-2,n}_{\infty} \cong H^{n+2k-2}_c(X_U)
\end{equation}

(Case 1) If $n$ is odd : Then $ E^{2k-2,n}_{\infty}=0$ also and hence we have $H^{n+2k-2}_c(X_U)=0$. So (\ref{eq:7}) gives an injection
$$
0 \to H^{n+2k-2}(X) \to  H^{n+2k-4}(Y)(-1)
$$

(Case 2) If $n$ is even : In this case, (\ref{eq:short}) gives
\begin{eqnarray*}
H^{n+2k-2}_c (X_U) &\cong& E^{2k-2,n}_{\infty}=E^{2k-2,n}_2 = H^{2k-2}_c(U,R^{n}{p_1}_* \Q) \\
&\cong& H^{2k-2}_c(U) \otimes H^n(Q_t) \cong H^{2k-2}(\PP^{k-1}) \otimes H^n(\PP^{n+k})
\end{eqnarray*}
by Lefschetz theorem. Hence, we can rewrite the exact sequence (\ref{eq:7}) as follows:
$$
0 \to H^{2k-2}(\PP^{k-1}) \otimes H^{n}(\PP^{n+k})\stackrel{h_*}{\to} H^{n+2k-2}(X) \to H^{n+2k-4}(Y) \to \cdots
$$
Recall that $H^{n+2k-2}_0(X) \cong H^{n+2k-2}(X)/\im~i^*_X$ and note that $\im ~h_* \cap H^{n+2k-2}_0(X) = \emptyset$: In fact, $h_*$ can be factorized as 
$$
\xymatrix{
H^{2k-2}(\PP^{k-1}) \otimes H^{n}(\PP^{n-k}) \ar@{^(->}[r] \ar@/^1pc/[rr]^{h_*}&H^{n+2k-2}(\PP^{k-1} \times \PP^{n+k}) \ar@{^(->}[r]_{\qquad i^*_X} & H^{n+2k-2}(X)
}
$$
Hence $\im~ h_* \subseteq \im ~i^*_X$, and we get an injection
$$
  0 \to H^{n+2k-2}_0(X) \to H^{n+2k-4}(Y)(-1) 
$$
 \end{proof}

Now we have a family $\pi:Y \to \Delta$ of quadrics of dimension $n+k-2$ which is even, so we can form a double covering of $\Delta$ as in the construction of O'Grady(\cite{ogrady}). For a general $t \in \Delta$, the fiber $\tilde{Q}_t$, which is a smooth quadric of dimension $n+k-2=2(m-1)$, contains two $\frac{m(m-1)}{2}-$dimensional irreducible families $F^1_t, F^2_t$ of $(m-1)-$planes\cite{gh}. As in the case of $n+k$ odd, we can form the following diagram
$$
\xymatrix@R=10pt{
&\ar[dl]_{pr_1} \Gamma \ar[dr]^{pr_2}& &\\
P \ar[d]_f \ar[drr]^{\psi} && Y \ar[d]^{\pi} \ar[r]^{\pi_2 \quad} & \PP^{n+k-1}\\
W \ar[rr]_{\sigma} && \Delta&
}
$$
where 
\begin{enumerate}
\item $P = \{ M \subset Y ~|~ \pi_1(M) = t \in \Delta \text{ a point}, \pi_2(M)\cong \PP^{m-1} \subset \PP^{n+k-1} \}$
\item $\Gamma = \{ (M, y) \in P \times Y~|~ \pi(y) =\pi(M) \in \Delta, \pi(y) \in \pi_2(M) \cong \PP^{m-1} \subset \tilde{Q}_{\pi(M)} \}$, $pr_1$ and $pr_2$ are projections.
\item $\psi : P \to \Delta$ is a natural map defined by $\psi(M) = \pi(M) =t \in \Delta$
\item $f: P \to W$ is $F-$bundle and $pr_1:\Gamma \to P$ is $\PP^{m-1}-$bundle, where $F$ is the abstract variety such that $F^i_t$ is isomorphic to it for $i=1,2$ and $t \in \Delta$
\item $\sigma: W \to \Delta$ is a double covering branched over the discriminant variety $\Delta_1$ of the family $\pi : Y\to \Delta$.
\end{enumerate}

\begin{lemma}\label{lemma:surjoddk}
There are surjections of Hodge structures
\begin{enumerate}
\item  $$ H^{k-2}(\tilde{W})(-m+1) \to H^{n+2k-4}(Y)  \qquad \text{if $n$ is odd}$$
where $\tilde{W}$ is a disjoint union of finitely many copies of $W$,
\item $$ \bigoplus_r H^{2r}(W)(-q_r)^{\oplus l_r} \to H^{n+2k-4}(Y) \qquad \text{if $n$ is even} $$
where $q_r=\frac{n+2k-4-2r}{2}$ and $ l_r=\sum_j \dim H_{n+2k-4-2j-2r}(F)$.
\end{enumerate}
\end{lemma}

\begin{proof}
First note that by using the same arguments in the proof of theorem \ref{thm:main} in the case of $n+k$ odd, we can show that 
\begin{equation}\label{eq:pr2}
{pr_2}_*: H_{n+2k-4}(\Gamma) \to H_{n+2k-4}(Y)
\end{equation}
is surjective and 
 \begin{eqnarray*}
 H_{n+2k-4}(\Gamma) &=& \bigoplus^{m-1}_{j=0} (H_{n+2k-4-2j} (P) \otimes H_{2j}(\PP^{m-1}) )\\
&=& \bigoplus^{m-1}_{j=0} H_{n+2k-4-2j} (P)(j)\\
&=& \bigoplus^{m-1}_{j=0} \bigoplus_{s}(H_{n+2k-4-2j-2s}(W) \otimes H_{2s}(F))(j) 
 \end{eqnarray*}

(Case 1) If $n$ is odd : In this case, $n+2k-4-2j-2s$ is odd. But $H_r(W)$ is zero for an odd number $r$ unless $r=k-2 =\dim W$. Hence we have
$$
H_{n+2k-4-2j}(P) =H_{k-2}(W) \otimes H_{n+k-2-2j}(F) = H_{k-2}(W)(m-1-j)^{\oplus l_j}
$$
where $l_j=\dim H_{n+k-2-2j}(F)$. Thus we have
\begin{eqnarray*}
  H_{n+2k-4}(\Gamma) 
&=& \bigoplus^{m-1}_{j=0} H_{k-2}(W)(m-1-j)^{\oplus l_j}(j)\\
&=&\bigoplus^{m-1}_{j=0} H_{k-2}(W)(m-1)^{\oplus l_j} = H_{k-2}(W)(m-1)^{\oplus N}
\end{eqnarray*}
where $N = \sum^m_{j=0} l_j$. Hence by (\ref{eq:pr2}) and using Poincar\'{e} duality, we get a surjection
$$
\xymatrix{
{pr_2}_* : H^{k-2}(W)(-m+1)^{\oplus N} \ar@{>>}[r] & H^{n+k-4}(Y)
}
$$
By setting $\tilde{W}$ to be a disjoint union of $N$ copies of $W$, we get the desired surjection in this case.

\medskip

(Case 2) If $n$ is even : In this case, $n+2k-4-2j$ is even. Hence we have
\begin{eqnarray*}
  H_{n+2k-4-2j}(P)&=& \bigoplus_{s}(H_{2s}(W) \otimes H_{n+2k-4-2j-2s}(F))
\end{eqnarray*}
Thus we have
\begin{eqnarray*}
  H_{n+2k-4}(\Gamma) &=& \bigoplus^{m-1}_{j=0} H_{n+2k-4-2j} (P)(j)\\
&=& \bigoplus^{m-1}_{j=0}\left(\bigoplus_{s}(H_{2s}(W) \otimes H_{n+2k-4-2j-2s}(F)) \right) (j)\\
&=& \bigoplus_{s} H_{2s}(W) \otimes \left(\bigoplus^{m-1}_{j=0}H_{n+2k-4-2j-2s}(F) \right)(j)\\
&=& \bigoplus_{s} H_{2s}(W)(q_s)^{\oplus l_s}
\end{eqnarray*}
where $q_s = \frac{n+2k-4-2s}{2}$ and $l_s = \sum^{m-1}_{j=0} \dim H_{n+2k-4-2j-2s}(F)$. Hence by (\ref{eq:pr2}) and using Poincar\'{e} duality, we get a surjection
$$
\xymatrix{
{pr_2}_* : \bigoplus_{r}H^{2r}(W)(-q_r)^{\oplus l_r}  \ar@{>>}[r]  & H^{n+2k-4}(Y)
}
$$
in this case. 
\end{proof}

We can summarize lemmas \ref{lemma:inj} and \ref{lemma:surjoddk} in the following diagram:
\begin{align}\label{diag:theta}
\xymatrix{
& H^* (W) \ar@{>>}[d]^{\theta} \ar[dr]^{\psi_2 \circ \theta}&\\
0 \to H^{n+2k-2}_0 (X) \ar[r]^{\psi_1} & H^{n+2k-4}(Y)(-1) \ar[r]^{\psi_2 \qquad \quad } & 
Gr^{W}_{n+2k-2}H^{n+2k-1}_c(X_U) \to 0 
}
\end{align}
where with the notation in theorem \ref{lemma:surjoddk}
$$
H^*(W) = \begin{cases} H^{k-2}(\tilde{W})(-m) & \text{if $n$ is odd} \\ 
\bigoplus_r H^{2r}(W)(-q_r-1)^{\oplus l_r} & \text{if $n$ is even} \end{cases}
$$
and $\theta: H^*(W) \to H^{n+2k-4}(Y)(-1)$ is the surjection in lemma \ref{lemma:surjoddk}.

\medskip

To finish the proof of theorem \ref{thm:main} in this case, we need to lift the surjection $\theta$ to a surjection onto $H^{n+2k-2}_0(X)$. In order to do that,
we construct a section $s_U$ of $\psi_2$ in geometric way. Following two lemmas will lead us the desired section $s_U$.

\begin{lemma}
  There is a decomposition
$$
H^{n+2k-4}(Y) \cong \begin{cases} H^{k-2}(\Delta, R^{n+k-2}\pi_*\Q) & \text{if $n$ is odd}\\
H^{k-2}(\Delta, R^{n+k-2}\pi_*\Q) \oplus \bigoplus^{k-2}_{p=0, p\neq \frac{k-2}{2}} (\alpha_p \otimes \beta_p)\Q & \text{if $n$ is even}
\end{cases}
$$ 
where $\alpha_p$ (resp. $\beta_p$) is an algebraic cycle generating $H^{2p}(\Delta)$ (resp. $H^{n+2k-4-2p}(\tilde{Q}_t)$).
\end{lemma}

\begin{proof}
 Consider the Leray spectral sequence associated to the map $\pi: Y \to \Delta$.
$$
E^{pq}_2 = H^p(\Delta, R^q {\pi}_* \Q) \quad \Rightarrow \quad H^{p+q}(Y, \Q)
$$
Note that 
\begin{equation}  \label{eq:ss}
  E^{pq}_2 = \begin{cases} 0 & \text{if either $q$ is odd or $p$ is odd, $p=k-2$} \\
H^p(\Delta,\Q) & \text{if $p$ is even and $q$ is even, $q \neq n+k-2$}\\
H^p (\Delta,R^{n+k-2}{\pi}_*\Q) &\text{if $p$ is even and $q=n+k-2$}  \end{cases}
\end{equation}
where $(R^{n+k}{\pi}_*\Q)_t=\Q^2$ for $t\in \Delta$. Then the spectral sequence degenerates at $E_2$. Let $F^{\dt}$ be the filtration on $H^{n+2k-4}(Y)$ obtained from the spectral sequence. Then for each $p=0,1,..., 2(k-2)$ we have a short exact sequence
\begin{equation}\label{eq:split}
0 \to F^{p+1} \to F^p \to \text{Gr}^p_F H^{n+2k-4}(Y) \to 0 
\end{equation}
where $$\text{Gr}^p_F H^{n+2k-4}(Y) = E^{p, n+2k-4-p}_{\infty} = H^{p}(\Delta , R^{n+2k-4-p}{\pi}_*\Q)$$

\medskip
(Case 1) If $n$ is odd (and hence $k$ is odd): then $n+2k-4$ is odd. Hence $E^{pq}_2=0$ unless $(p,q)=(k-2, n+k-2)$. Hence we have
\begin{equation}\label{eq:xdelta}
H^{n+2k-4}(Y) \cong H^{k-2}(\Delta, R^{n+k-2}{\pi}_*\Q)
\end{equation}

\medskip
(Case 2) If $n$ is even (and hence $k$ is even): then $n+2k-4$ is also even. Note $ F^{2p-1} = F^{2p}$ from (\ref{eq:ss}). Hence we may write above sequence as
\begin{equation}\label{seq:splitting}
0 \to F^{2p+2} \to F^{2p} \stackrel{t_{2p}}{\to} \text{Gr}^{2p}_F H^{n+2k-4}(Y) \to 0
\end{equation}
First by using descending induction on $p$, we show that for $p$ such that $k-2 < 2p \leq 2(k-2)$ there is a natural splitting of exact sequences (\ref{seq:splitting}) such that 
$$
F^{2p} \cong F^{2p+2} \oplus \text{Gr}^{2p}_FH^{n+2k-4}(Y) \cong \bigoplus^{k-2}_{l=p+1} (\alpha_l \otimes \beta_l) \Q
$$
where $\alpha_l$ is an algebraic cycle generating $H^{2l}(\Delta) \cong \Q$ and $\beta_l$ is an algebraic cycle generating $ H^{n+2k-4-2l}(\tilde{Q}_t)(-1) \cong \Q$. 

If $2p=2(k-2)$, then 
\begin{eqnarray*}
F^{2(k-2)} &\cong& \text{Gr}^{2(k-2)}_F H^{n+2k-4}(Y) = E^{2(k-2), n}_{\infty}\\
&=& H^{2(k-2)}(\Delta) \otimes H^{n}(\tilde{Q}_t) \cong (\alpha_{k-2} \otimes \beta_{k-2})\Q
\end{eqnarray*}
where $\alpha_{k-2}$(resp. $\beta_{k-2}$) is an algebraic cycle generating $H^{2(k-2)}(\Delta)=\Q$ (resp. $ H^{n}(\tilde{Q}_t)=\Q$). Now suppose $k-2<2p<2(k-2)$ and consider the short exact sequnce (\ref{seq:splitting}). Since $2p>k-2$, $\dim H^{2p}(\Delta)=\dim H^{n+2k-4-2p}(\tilde{Q}_t)=1$, hence we can choose an algebraic cycle $\alpha_p$ (resp. $\beta_p$) which generates $H^{2p}(\Delta)$ (resp. $H^{n+2k-4-2p}(\tilde{Q}_t)$) and hence
\begin{eqnarray*}
  \text{Gr}^{2p}_F H^{n+2k-4}(Y) &=&H^{2p}(\Delta, R^{n+2k-4-2p}{\pi}_*\Q) = H^{2p}(\Delta) \otimes H^{n+2k-4-2p}(\tilde{Q}_t) \\
&=& (\alpha_p \otimes \beta_p)\Q
\end{eqnarray*}
and by induction hypothesis $F^{2p+2}$ has a decomposition
$$
F^{2p+2} \cong \bigoplus^{k-2}_{l=p+1}(\alpha_l \otimes \beta_l) \Q
$$
with a basis $\mathcal{B}_{p+1}=\{\alpha_{p+1} \otimes \beta_{p+1},..., \alpha_{k-2} \otimes \beta_{k-2}\}$, where $a_l$ and $\beta_l$ are algebraic cycles. Let $\mathcal{B}'_{p}=\{v_{p}, \alpha_{p+1} \otimes \beta_{p+1},..., \alpha_{k-2} \otimes \beta_{k-2}\}$ be a basis of $F^{2p}$ extending the basis $\mathcal{B}_{p+1}$. Then $t_{2p} (v_p ) = q_p (\alpha_p \otimes \beta_p)$ for some nonzero $q_p \in \Q$. Define 
$$
s_p : \text{Gr}^{2p}_F H^{n+2k-4}(Y) \to F^{2p}
$$
by $s_p(\alpha_p \otimes \beta_p) = \frac{1}{q_p}v_p$. Then $s_p$ is a section of $t_p$ and hence we have a decomposition
$$
F^{2p}\cong F^{2p+2} \oplus \text{Gr}^{2p}_F H^{n+2k-4}(Y) =\bigoplus^{k-2}_{l=p+1}(\alpha_l \otimes \beta_l)\Q \oplus s_p(\alpha_p \otimes \beta_p)\Q
$$
with a basis $\mathcal{B}_p=\{\alpha_p \otimes \beta_p, \alpha_{p+1} \otimes \beta_{p+1}, ..., \alpha_{k-2}\otimes \beta_{k-2} \}$ of $F^{2p}$ by identifying $\alpha_p \otimes \beta_p$ with its image under $s_p$.

Now for splitting for $2p < k-2$ in the exact sequence (\ref{seq:splitting}), the argument is same as the case when $2p>k-2$, since $\text{Gr}^{2p}_F H^{n+2k-4}(Y) = H^{2p}(\Delta) \otimes H^{n+2k-4-2p}(\tilde{Q}_t) = \Q$. Hence we can choose a natural splitting of each exact sequence (\ref{seq:splitting}) and get a decomposition
\begin{equation}\label{eq:decomp}
H^{n+2k-4}(Y) \cong H^{k-2}(\Delta,R^{n+k-2}{\pi}_*\Q) \oplus \bigoplus^{k-2}_{p=0, p\neq \frac{k-2}{2}} (\alpha_p \otimes \beta_p)\Q
\end{equation}
\end{proof}

\begin{lemma}

Let $L$ be the subspace of $H^{n+2k-4}(Y)$ generated by $\{\alpha_p \otimes \beta_p ~|~ p =0,...,\widehat{\frac{k-2}{2}},..., k-2\}$. Then,
  $$
H^{k-2}(\Delta, R^{n+k-2}\pi_* \Q) \cong L^{\perp}
$$
where $\perp$ is the orthogonal complement with respect to the cup product on $H^{n+2k-4}(Y)$.
\end{lemma}
\begin{proof}
First we show that
$$
(\alpha_i \otimes \beta_i) \cup (\alpha_j \otimes \beta_j) = 
\begin{cases} 0 & \text{if $i+j \neq k-2$} \\ \deg~ \Delta & \text{if $i+j=k-2$}\end{cases}
$$  
Since
\begin{multline}
H^{2i}(\Delta, R^{n+2k-4-2i}\pi_*\Q) \otimes H^{2j}(\Delta, R^{n+2k-4-2j}\pi_*\Q)\\
 \stackrel{\cup}{\longrightarrow} \quad H^{2(i+j)}(\Delta, R^{2n+4k-8-2(i+j)}\pi_*\Q)
\end{multline}
$(\alpha_i \otimes \beta_i) \cup (\alpha_j \otimes \beta_j) =0$ for $i+j \neq k-2$. 
For $i+j=k-2$, we observe $\alpha_i, \beta_i$ closely. Note that for any $i \neq \frac{k-2}{2}$, $\alpha_i = h^i_{\Delta}$ where $h^i_{\Delta} \in H^{2i}(\Delta)$ is a class corresponding to an iterated hyperplane section of $\Delta$ of codimension $i$. For $\beta_i$, let  $\Pi \in H^{n+k-2}(\tilde{Q}_t)$ be a class corresponding to $\PP^{\frac{n+k-2}{2}} \subset \tilde{Q}_t$ and  $H^r \in H^{2r}(\tilde{Q}_t)$ a class corresponding to an iterated hyperplane section of $\tilde{Q}_t$ of codimension $r$. Note that $\Pi$ can be chosen in either families $F^1_t$ or $F^2_t$ of $\PP^{\frac{n+k-2}{2}}$ contained in $\tilde{Q_t}$ (\cite{reid}). Then, 
$$
\beta_i = \begin{cases}  H^{\frac{n+2k-4-2i}{2}} & \text{if $i > \frac{k-2}{2}$} \\
\Pi \cup H^{\frac{k-2}{2}-i} & \text{if $i<\frac{k-2}{2}$}
\end{cases}
$$
Then for $i+j=k+2$, 
$$
(\alpha_i \otimes \beta_i ) \cup (\alpha_j \otimes \beta_j)=\deg ~\Delta
$$

Now we show $L \cap L^{\perp}=0$. Let $\eta=\sum^{k-2}_{i=0,i\neq \frac{k-2}{2}} c_i (\alpha_i \otimes \beta_i) \in L \cap L^{\perp}$. Then, we have $c_i \cdot \text{deg} \Delta = 0 $ for any $i \neq \frac{k-2}{2}$. Hence $\eta =0$. Hence
$$
H^{n+2k-4}(Y) = L \oplus H^{k-2}(\Delta, R^{n+k-2}\pi_* \Q) = L \oplus L^{\perp}
$$
Hence we get
$$
H^{k-2}(\Delta, R^{n+k-2}\pi_* \Q )= L^{\perp}
$$
\end{proof}

From the above lemmas, we have an injection 
\begin{equation}\label{eq:h}
h: H^{k-2} (\Delta, R^{n+k} {p_1}_*\Q) \hookrightarrow H^{n+2k-2}(X_{\Delta})
\end{equation}

We refer the diagram (\ref{diag:theta}) for the following lemma.

\begin{lemma}\label{lemma:section}
 There is a section defined in geometric way
$$s_U:  Gr^{W}_{n+2k-2}H^{n+2k-1}_c(X_U) \to H^{n+2k-4}(Y)(-1)$$ of $\psi_2$.
\end{lemma}

\begin{proof}
First recall $H^{n+2k-4}(Y)(-1) \cong H^{n+2k-2}(X_{\Delta})$ and $\psi_2$ is the connecting homomorphism in the exact sequence of mixed Hodge structures (\ref{eq:mhs}) 
$$
\cdots \to H^{n+2k-2}(X) \to H^{n+2k-2}(X_{\Delta}) \stackrel{\psi_2}{\to} Gr^W_{n+2k-2} H^{n+2k-1}_c(X_U) \to 0
$$
As in the proof of lemma \ref{lemma:inj}, we can show that the the Leray spectral sequence associated to ${p_1} : X_U \to U$ gives rise to an exact sequence
$$
0 \to E^{2k-2,n+1}_{\infty} \to H^{n+2k-1}_c(X_U) \to E^{k-1,n+k}_{\infty} \to 0
$$ 
By taking the exact functor $\text{Gr}^W_{n+2k-2}$, we get  an isomorphism 
$$
l:Gr^W_{n+2k-2}H^{n+2k-1}(X_U) \stackrel{\cong}{\longrightarrow} Gr^W_{k-2}H^{k-1}_c(U)
$$
since $E^{2k-2,n+1}_{\infty}=H^{2k-2}_c(U)(-\frac{n+1}{2}) \cong H^{2k-2}(\PP^{k-1})(-\frac{n+1}{2})$ is a pure Hodge structure of weight $n+2k-2,$ if $n$ is odd and $0$ if $n$ is even. Hence we have a following commutative diagram:
\begin{align}\label{diag:section}
\xymatrix{
H^{n+2k-2}(X_{\Delta}) \ar[r]^{\psi_2 \qquad }  & Gr^W_{n+2k-2} H^{n+2k-1}_c(X_U) \ar[d]^{\cong}_l\\
H^{k-2}(\Delta, R^{n+k}{p_1}_* \Q) \ar@{^(->}[u]^{h} 
\ar[r]^{\delta \quad } & Gr^W_{k-2}H^{k-1}_c(U) \left(-\frac{n+k}{2} \right)
}
\end{align}
where $h$ is the injection of (\ref{eq:h}) and $\delta: H^{k-2}(\Delta, R^{n+k}{p_1}_* \Q)\to Gr^W_{k-2}H^{k-1}_c(U) \left(-\frac{n+k}{2} \right)$ is the connecting homomorphism of long exact sequence induced by a short exact sequence of sheaves on $\PP^{k-1}$
$$
0 \to j_! (R^{n+k}p_1 \Q)|_U \to R^{n+k}{p_1}_*\Q \to {i_{\Delta}}_*(R^{n+k}{p_1}_* \Q)|_{\Delta} \to 0
$$
where $j: U \to \PP^{k-1}$ and $i_{\Delta} : \Delta \to \PP^{k-1}$ are inclusions. Hence in order to choose a section of $\psi_2$, it is enough to construct a section of $\delta$. 

Note that $(R^{n+k}{p_1}_*\Q )_t =\Q^2$ for $t\in \Delta$ and $R^{n+k}{p_1}_*\Q|_U =\Q_U$ is the constant sheaf. Let $G$ be the monodromy group of $R^{n+k}{p_1}_*\Q|_{\Delta}$. Then $R^{n+k}{p_1}_*\Q|^G_{\Delta} =\Q$ is a subsheaf of $R^{n+k}{p_1}_*\Q|_{\Delta}$ and let $R^{n+k}{p_1}_*\Q|^G_{\Delta} \hookrightarrow R^{n+k}{p_1}_*\Q|_{\Delta}$ be an injection. Let 
$$
s: {i_{\Delta}}_* (R^{n+k}{p_1}_* \Q^G_{\Delta}) \to {i_{\Delta}}_* (R^{n+k}{p_1}_* \Q_{\Delta}) 
$$ 
be the induced map. Then there is a sheaf $s^*(R^{n+k}{p_1}_*\Q)$ on $\PP^{k-1}$ which fits into the following commutative diagram:
$$
\xymatrix{
0 \ar[r]&   j_! (R^{n+k}p_1 \Q|_U) \ar[r] \ar@{=}[d] & s^*(R^{n+k}{p_1}_*\Q) \ar[d] \ar[r] & {i_{\Delta}}_*(R^{n+k}{p_1}_* \Q |_{\Delta}^{G}) \ar[d]^{s} \to 0 \\
0 \ar[r] &  j_! (R^{n+k}p_1 \Q |_U) \ar[r]& R^{n+k}{p_1}_*\Q \ar[r]& {i_{\Delta}}_*(R^{n+k}{p_1}_* \Q|_{\Delta}) \to 0 \\
}
$$
From this, we get a commutative diagram
$$
\xymatrix@C=8pt{
\cdots H^{k-2}(\PP^{k-1},s^*(R^{n+k}{p_1}_*\Q)) \ar[d] \ar[r] &H^{k-2}(\Delta, \Q) \ar[r]^{\delta_1\qquad } \ar[d]^{s} & Gr^W_{k-2}H^{k-1}_c(U,\Q)  \ar@{=}[d] \to 0 \\
\cdots H^{k-2}(\PP^{k-1},R^{n+k}{p_1}_*\Q)\ar[r] &H^{k-2}(\Delta, R^{n+k}{(p_1|_{\Delta})}_* \Q) \ar[r]^{\quad \delta} & Gr^W_{k-2}H^{k-1}_c(U) \to 0 
}
$$
i.e.
\begin{equation}\label{eq:delta}
\delta_1 = \delta \circ s
\end{equation}
Now from the exact sequence
$$
\cdots \to H^{k-2}(\PP^{k-2}) \stackrel{i^*_{\Delta}}{\to} H^{k-2}(\Delta) \stackrel{\delta_1}{\to} Gr^{W}_{k-2}H^{k-1}_c(U) \to 0
$$
we have an isomorphism
$$
Gr^{W}_{k-2}H^{k-1}_c(U) \stackrel{\cong}{\longrightarrow} H^{k-2}(\Delta)/\im ~ i^*_{\Delta}
$$
hence we can choose a natural section $s_1 : Gr^W_{k-2}H^{k-1}_c(U) \to H^{k-2}(\Delta)$ of $\delta_1 : H^{k-2}(\Delta) \to Gr^W_{k-2}H^{k-1}_c(U)$, i.e. $\delta_1 \circ s_1 =\text{id}$. Then from (\ref{eq:delta}) we get a section $s \circ s_1$ of $\delta$. By combining all these, we get a commutative diagram:

$$
\xymatrix@R=13pt{
H^{n+2k-2}(X_{\Delta}) \ar[r]^{\psi_2 \qquad }  & Gr^W_{n+2k-2} H^{n+2k-1}_c(X_U) \ar[d]^{\cong}_l\\
H^{k-2}(\Delta, R^{n+k}{p_1}_* \Q) \ar@{^(->}[u]^{h} 
\ar[r]^{\delta \quad } & Gr^W_{k-2}H^{k-1}_c(U) \left(-\frac{n+k}{2} \right) 
\\
H^{k-2}(\Delta) \ar[r]^{\delta_1 \qquad } \ar[u]^{s}& \ar@/^1pc/[l]^{s_1} Gr^W_{k-2}H^{k-1}_c(U) \left(-\frac{n+k}{2} \right) \ar@{=}[u]
}
$$
Set 
$$
s_U = h \circ s \circ s_1 \circ l
$$ 
then 
$$
\psi_2 \circ s_U = \psi_2 \circ( h \circ s \circ s_1 \circ l) = l^{-1} \circ (\delta \circ s \circ s_1) \circ l = l^{-1}\circ l =\text{id} 
$$
Hence $s_U$ is a section of $\psi_2$.
\end{proof}

\begin{proof}[Proof of theorem \ref{thm:main} when $n+k$ is even]
  Now we finish the proof of the Theorem \ref{thm:main} in case when $n+k$ is even. Recall the diagram (\ref{diag:theta})
$$
\xymatrix{
& 
H^* (W) \ar@{>>}[d]^{\theta} \ar[dr]^{\psi_2 \circ \theta}&\\
0 \to H^{n+2k-2}_0 (X) \ar[r]^{\psi_1} & H^{n+2k-4}(Y)(-1) \ar[r]^{\psi_2 \qquad \quad } & \ar@/^ 1pc/[l]^{s_U} Gr^{W}_{n+2k-2}H^{n+2k-1}_c(X_U) \to \cdots
}
$$
where $s_U$ is the section of $\psi_2$ constructed in lemma \ref{lemma:section}.

Let 
$$\Theta = \theta - s_U \circ \psi_2 \circ \theta : H^{*}(W) \to H^{n+2k-4}(Y)(-1)$$
First note that $ \psi_2 \circ \Theta =0 $: In fact,
\begin{eqnarray*}
\psi_2 \circ \Theta &=& \psi_2 \circ (\theta - s_U \circ \psi_2 \circ \theta) \\
&=& \psi_2 \circ \theta - \psi_2 \circ s_U \circ \psi_2 \circ \theta =\psi_2 \circ \theta - \psi_2 \circ \theta=0
\end{eqnarray*}
Hence $\im ~\Theta \subseteq \ker ~\psi_2$ and we may consider that $\Theta$ is mapped into $H^{n+2k-2}_0(X)$ since $\psi_1$ is an injection. 
$$
\xymatrix{
& \ar[dl]_{\Theta=\theta - s_U \circ \psi_2 \circ \theta} H^* (W) \ar@{>>}[d]^{\theta} \ar[dr]^{\psi_2 \circ \theta}&\\
0 \to H^{n+2k-2}_0 (X) \ar[r]^{\psi_1} & H^{n+2k-4}(Y)(-1) \ar[r]^{\psi_2 \qquad \quad } & \ar@/^ 1pc/[l]^{s_U} Gr^{W}_{n+2k-2}H^{n+2k-1}_c(X_U) \to \cdots
}
$$
Now we show this map is in fact a surjection. Let $\alpha \in H^{n+2k-2}_0(X)$. Since $\psi_1$ is an injection we may identify $\alpha$ with $\psi_1(\alpha)$. Then there is $\beta \in H^*(W)$ such that $\theta(\beta) = \psi_1(\alpha)$. Then
\begin{eqnarray*}
\Theta (\beta) &=& (\theta-s_U \circ \psi_2 \circ \theta)(\beta) \\
&=& \theta(\beta) -s_U \circ \psi_2 \circ \theta(\beta)=\psi_1(\alpha)-s_U \circ \psi_2 \circ \psi_1(\alpha) = \psi_1(\alpha) =\alpha
\end{eqnarray*}
Hence we have a surjection 
$$
\xymatrix{
\Theta : H^*(W) \ar@{>>}[r] & H^{n+2k-2}_0(X)
}
$$
as we claimed. This finishes the proof of theorem \ref{thm:main}.
\end{proof}

\end{document}